\documentclass[reqno,12pt]{amsart}
\usepackage{amsmath,amsthm,amssymb,amsfonts,amscd,graphicx}
\theoremstyle{plain} 
	\newtheorem{thm}{Theorem}[section]
	\newtheorem*{thm*}{Theorem}
	\newtheorem{cor}[thm]{Corollary}
	\newtheorem*{cor*}{Corollary}

	\newtheorem{prop}[thm]{Proposition}
	\newtheorem*{prop*}{Proposition}
	\newtheorem{conj}[thm]{Conjecture}
	\newtheorem*{conj*}{Conjecture}
	
\theoremstyle{definition}
	\newtheorem{defn}[thm]{Definition}

\theoremstyle{remark}
	\newtheorem{rem}[thm]{Remark}
	
	\newtheorem*{pf}{Proof}
\numberwithin{equation}{section}

\def\CC{{\mathbb C}}

\def\PP{{\mathbb P}}
\def\QQ{{\mathbb Q}}
\def\RR{{\mathbb R}}
\def\ZZ{{\mathbb Z}}

\def\B{{\mathcal B}}
\def\C{{\mathcal C}}
\def\D{{\mathcal D}}
\def\E{{\mathcal E}}
\def\F{{\mathcal F}}

\def\H{{\mathcal H}}
\def\I{{\mathcal I}}

\def\K{{\mathcal K}}
\def\L{{\mathcal L}}
\def\M{{\mathcal M}}

\def\O{{\mathcal O}}

\def\U{{\mathcal U}}
\def\V{{\mathcal V}}

\def\T{{\mathcal T}}

\def\m{{\mathfrak m}}
\def\p{\partial }

\def\ns{{\nabla}\hspace{-1.4mm}\raisebox{0.3mm}{\text{\footnotesize{\bf /}}}}

\begin{document}
\title[Stokes matrices for orbifold projective lines]
{Stokes Matrices for the Quantum Cohomologies of Orbifold Projective Lines}
\date{\today}
\author{Kohei Iwaki}
\address{Research Institute for Mathematical Sciences, Kyoto University, 
Kyoto, 606-8502, Japan}
\email{iwaki@kurims.kyoto-u.ac.jp}
\author{Atsushi Takahashi}
\address{Department of Mathematics, Graduate School of Science, Osaka University, 
Toyonaka Osaka, 560-0043, Japan}
\email{takahashi@math.sci.osaka-u.ac.jp}
\begin{abstract}
We prove the Dubrovin's conjecture for the Stokes matrices for the quantum cohomology 
of orbifold projective lines. The conjecture states that the Stokes matrix of the 
first structure connection of the Frobenius manifold constructed from the 
Gromov-Witten theory coincides with the Euler matrix of a full exceptional collection 
of the bounded derived category of the coherent sheaves. Our proof is based on 
the homological mirror symmetry, primitive forms of affine cusp polynomials 
and the Picard-Lefschetz theory.
\end{abstract}
\maketitle
\section{Introduction}

For a smooth projective variety $X$ over $\CC$, the quantum cohomology ring of $X$ 
is defined as a generalization of usual cohomology ring. The quantum cohomology ring 
coincide with the usual cohomology ring as a vector space, but the product structure is 
``quantum corrected'' from the usual cup product, by counting the number of 
holomorphic curves in $X$ hitting the cycles Poincar\'e dual to the cohomology classes. 
Such counting numbers are called the Gromov-Witten invariants of $X$. This idea comes 
from physics and attracted the interests of mathematicians because of spectacular 
predictions for classical enumeration problems in algebraic geometry through 
the mirror symmetry. Now the theory of quantum cohomology and the Gromov-Witten 
invariants are extended to cases when $X$ is an orbifold \cite{agv:1, cr:1}.

The theory of Frobenius manifold formulated by Dubrovin \cite{d:1}, which first 
appeared as the flat structure in K.~Saito's study of the deformation space of 
an isolated hypersurface singularity, 
enable us to treat quantum cohomology systematically. A Frobenius manifold is 
a complex (formal) manifold whose tangent space at any point has a bilinear form 
and an associative commutative product with certain compatibility conditions.  
From these compatibility conditions, it turns out to be that, the structure constants 
of the product structure are given by third derivatives of a function on the 
Frobenius manifold. The function is called Frobebius potential. The quantum cohomology 
of $X$ satisfies the axioms of Frobenius manifolds, and its Frobenius potential 
is the generating function of the genus zero Gromov-Witten invariants of $X$.

The above compatibility conditions give a property of integrable systems 
to Frobenius manifolds. Namely, for any Frobenius manifold $M$, we can construct 
a flat connection on the tangent bundle of $M\times\PP^{1}$, called the
{\it first structure connection} of $M$. The differential equation satisfied 
by the flat sections of the first structure connection can be regarded as 
an isomonodromic family of a meromorphic ordinary differential equation on $\PP^{1}$,
parametrized by points of $M$. The equation has two singular points on $\PP^{1}$; 
one is a regular singular point and the other one is an irregular singular point.  
If a point on $M$ is semi-simple, that is, if there are no nilpotent elements 
in the product structure on the tangent space at the point, then one can define 
the {\it monodromy data} of the differential equation at the point,
consisting of the monodromy matrix at the regular singular point, 
the Stokes matrix at the irregular singular point, and the connection matrix between 
these two singular points. These data do not depend on the choice of a semi-simple point 
due to the isomonodromy property. Moreover, as is shown in \cite{d:2}, from these data 
we can reconstruct the Frobenius structure by the Riemann-Hilbert correspondence. 
 
After the Zaslow's work \cite{z}, Dubrovin formulated a conjecture in 
\cite{d:3} about a close relationship between the monodormy data of the 
Frobenius manifold constructed from the Gromov-Witten theory of $X$ 
and the structure of the bounded derived category of coherent sheaves on $X$. 
To state his conjecture, we recall some notions here.
\begin{defn}\label{defn:exceptional}
Let $\T$ be a $\CC$-linear triangulated category $\T$ with a translation functor $T$.
\begin{enumerate}
\item
An object $\E$ in $\T$ is called an {\it exceptional object} 
(or is called {\it exceptional}) if $\T(\E,\E)=\CC\cdot {\rm id}_\E$ and 
$\T(\E,T^p\E)=0$ when $p\ne 0$.
\item 
An {\it exceptional collection} in $\T$ is an ordered set $(\E_1,\dots,\E_\mu)$ 
of exceptional objects satisfying the condition $\T(\E_i,T^p\E_j)=0$
for all $p\in\ZZ$ and $i>j$.
\item 
An exceptional collection $(\E_1,\dots,\E_\mu)$ in $\T$ is called {\it full}
if the smallest full triangulated subcategory of $\T$ containing $\E_1,\dots,\E_\mu$
is equivalent to $\T$.
\end{enumerate}
\end{defn}
\begin{defn}
Let $\T$ be a $\CC$-linear triangulated category $\T$ with a translation functor $T$.
Assume that $\T$ is finite, namely, for all objects $\E,\E'\in \T$ one has 
\begin{equation}
\sum_{p\in\ZZ}\dim_\CC\T(\E,T^p\E') <\infty.
\end{equation}
\begin{enumerate}
\item
Let $K_0(\T)$ be the Grothendieck group of $\T$. 
The pairing $\chi:K_0(\T)\times K_0(\T)\longrightarrow \ZZ$ defined by 
\begin{equation}
\chi\left(\left[\E\right],\left[\E'\right]\right):=\sum_{p\in\ZZ}(-1)^p\dim_\CC\T(\E,T^p\E')
\end{equation}  
is called the {\it Euler pairing}. 
\item
Suppose that  $\T$ is generated by a full exceptional collection 
$(\E_1,\dots,\E_\mu)$. We shall call the $\mu\times \mu$-matrix $\chi=(\chi_{ij})$, 
$\chi_{ij}=\chi\left(\left[\E_i\right],\left[\E_j\right]\right)$ 
the {\it Euler matrix} of the exceptional collection $(\E_1,\dots,\E_\mu)$.
\end{enumerate}
\end{defn}
Then, (a part of) the Dubrovin's conjecture is formulated as follows$:$
\begin{conj}[\cite{d:3}]
Let $X$ be a smooth projective variety over $\CC$ and let $M$ be the Frobenius manifold 
constructed from the Gromov--Witten theory for $X$
\begin{enumerate}
\item
The Frobenius manifold $M$ is semi-simple if and only if the 
bounded derived category $D^b{\rm coh}(X)$ of coherent sheaves on $X$ 
admits a full exceptional collection $(\E_1,\dots,\E_\mu)$ where 
$\mu:=\displaystyle\sum_{p\in \ZZ}\dim_\CC H^{p,p}(X)$. 
\item
Fix a semi-simple point $t\in M$ such that the values of the canonical coordinates at 
the point $t$ are pairwise distinct. Then the Stokes matrix of the first structure 
connection of the Frobenius structure defined at the point $t$ 
is identified with the Euler matrix of an exceptional collection 
$(\E_1,\dots,\E_\mu)$ in $D^b{\rm coh}(X)$.
\end{enumerate}
\end{conj}
The conjecture is proved for some examples of smooth projective varieties
(cf.~\cite{d:2, guz, u:1, u:2}).
In this paper, we shall consider a generalization of Dubrovin's conjecture 
to an orbifold projective line $\PP^1_A:=\PP^1_{(a_1,a_2,a_3)}$, 
an orbifold $\PP^1$ with at most three isotropic points of orders $a_1,a_2,a_3$ 
satisfying the condition $\chi_A:=1/a_1+1/a_2+1/a_3-1>0$. 
The existence of a full exceptional collection in $D^b{\rm coh}(\PP^1_A)$ 
is already shown by Geigle--Lenzing (Proposition~4.1 in \cite{gl:1}).
The semi-simplicity of the quantum cohomology ring of $\PP^1_A$ is a corollary 
(see Corollary~\ref{cor:semi-simplicity} below) of the classical mirror symmetry, 
the isomorphism of Frobenius manifolds between the one constructed from 
Gromov--Witten theory for $\PP^1_A$ and the one constructed from the theory of 
primitive forms for the polynomial $f_A:=x_1^{a_1}+x_2^{a_2}+x_3^{a_3}-x_1x_2x_3$.
The following is our main theorem$:$
\begin{thm}
Let $M_{\PP^1_A}$ be the Frobenius manifold constructed from the Gromov--Witten theory 
for $\PP^1_A$. Fix a semi-simple point $t\in M_{\PP^1_A}$ such that the values of 
the canonical coordinates at the point $t$ are pairwise distinct.
Then the Stokes matrix of the first structure connection of the Frobenius structure 
defined at the point $t$ is identified with the Euler matrix of a full exceptional 
collection in $D^b{\rm coh}(\PP^1_A)$.
\end{thm}

The proof is based on the homological mirror symmetry; an triangulated equivalence 
between $D^{b}{\rm coh}(\PP^{1}_{A})$ and $D^{b}{\rm Fuk}^{\rightarrow}(f_{A})$. 
The latter is the bounded derived category of the directed Fukaya category 
${\rm Fuk}^{\rightarrow}(f_{A})$ of $f_{A}$. A full exceptional collection of the 
$D^{b}{\rm Fuk}^{\rightarrow}(f_{A})$ is given by the vanishing cycles in the 
fiber of a suitable deformation of $f_{A}$, and the Euler forms of them can be written 
in terms of intersection numbers of these vanishing cycles multiplied by $-1$. 
On the other hand, due to the classical mirror symmetry, the oscillatory integral 
of the primitive form for $f_{A}$ along a certain cycle called Lefschetz thimble, 
gives a flat sections of the first structure connection of the Frobenius manifold 
constructed from $\PP^{1}_{A}$. Using the idea of \cite{u:2} and the Picard-Lefschetz 
theory, we can compute the Stokes matrix of the oscillatory integrals by the 
intersection numbers of vanishing cycles. Thus we obtain the main theorem. 
\bigskip
\noindent
{\it Acknowledgement}\\
\indent
The first named author thanks to Hiroshi Iritani and Yuuki Shiraishi for 
helphul discussions. The second named author is supported 
by JSPS KAKENHI Grant Number 24684005. 
\section{Preliminary}
\subsection{Definition of the Frobenius manifold}

In this section, we recall the definition and some basic properties of 
the Frobenius manifold \cite{d:1}. The definition below is taken 
from Saito-Takahashi \cite{st:1}.
\begin{defn}
Let $M=(M,\O_{M})$ be a connected complex manifold or 
a formal manifold over $\CC$ of dimension $\mu$
whose holomorphic tangent sheaf and cotangent sheaf 
are denoted by $\T_{M}, \Omega_M^1$ respectively
and $d$ be a complex number.
A {\it Frobenius structure of rank $\mu$ and dimension $d$ on M} is a tuple 
$(\eta, \circ , e,E)$, where $\eta$ is a non-degenerate $\O_{M}$-symmetric 
bilinear form on $\T_{M}$, $\circ $ is $\O_{M}$-bilinear product on $\T_{M}$, 
defining an associative and commutative $\O_{M}$-algebra structure with 
the unit $e$, and $E$ is a holomorphic vector field on $M$, called
the Euler vector field, which are subject to the following axioms:
\begin{enumerate}
\item The product $\circ$ is self-ajoint with respect to $\eta$: that is,
\begin{equation*}
\eta(\delta\circ\delta',\delta'')=\eta(\delta,\delta'\circ\delta''),\quad
\delta,\delta',\delta''\in\T_M. 
\end{equation*} 
\item The {\rm Levi}--{\rm Civita} connection $\ns:\T_M\otimes_{\O_M}\T_M\to\T_M$ 
with respect to $\eta$ is flat: that is, 
\begin{equation*}
[\ns_\delta,\ns_{\delta'}]=\ns_{[\delta,\delta']},\quad \delta,\delta'\in\T_M.
\end{equation*}
\item The tensor $C:\T_M\otimes_{\O_M}\T_M\to \T_M$  defined by 
$C_\delta\delta':=\delta\circ\delta'$, $(\delta,\delta'\in\T_M)$ is flat: that is,
\begin{equation*}
\ns C=0.
\end{equation*} 
\item The unit element $e$ of the $\circ $-algebra is a 
$\ns$-flat homolophic vector field: that is,
\begin{equation*}
\ns e=0.
\end{equation*} 
\item The metric $\eta$ and the product $\circ$ are homogeneous of degree 
$2-d$ ($d\in\CC$), $1$ respectively with respect to Lie derivative 
$Lie_{E}$ of {\rm Euler} vector field $E$: that is,
\begin{equation*}
Lie_E(\eta)=(2-d)\eta,\quad Lie_E(\circ)=\circ.
\end{equation*}
\end{enumerate}
A manifold $M$ equipped with a Frobenius structure $(\eta, \circ , e,E)$ is called a {\it Frobenius manifold}.
\end{defn}
From now on in this section, we shall always denote by $M$ a Frobenius manifold.
We expose some basic properties of Frobenius manifolds without their proofs.
Let us consider the space of horizontal sections of the connection $\ns$:
\[
\T_M^f:=\{\delta\in\T_M~|~\ns_{\delta'}\delta=0\text{ for all }\delta'\in\T_M\}
\]
which is a local system of rank $\mu $ on $M$ such that the metric $\eta$ 
takes constant value on $\T_M^f$. Namely, we have 
\begin{equation}
\eta (\delta,\delta')\in\CC,\quad  \delta,\delta' \in \T_M^f.
\end{equation}
\begin{prop}\label{prop:flat coordinates}
At each point of $M$, there exist a local coordinate $(t_1,\dots,t_{\mu})$, 
called flat coordinates, such that $e=\p_1$, $\T_M^f$ is spanned by 
$\p_1,\dots, \p_{\mu}$ and $\eta(\p_i,\p_j)\in\CC$ for all $i,j=1,\dots, \mu$,
where we denote $\p/\p t_i$ by $\p_i$. 
\end{prop} 
The axiom $\ns C=0$, implies the following:
\begin{prop}\label{prop:potential}
At each point of $M$, there exist the local holomorphic function $\F$, 
called Frobenius potential, satisfying
\begin{equation*}
\eta(\p_i\circ\p_j,\p_k)=\eta(\p_i,\p_j\circ\p_k)=\p_i\p_j\p_k \F,
\quad i,j,k=1,\dots,\mu,
\end{equation*}
for any system of flat coordinates. In particular, one has
\begin{equation*}
\eta_{ij}:=\eta(\p_i,\p_j)=\p_1\p_i\p_j \F. 
\end{equation*}
\end{prop}
The product structure on $\T_{M}$ is 
described locally by ${\mathcal F}$ as 
\begin{equation} \label{eq:prod}
{\partial_{i}} \circ 
{\partial_{j}}  = \sum_{k = 1}^{\mu}
c_{ij}^{k} {\partial_{k}}  
\quad i,j = 1, \cdots, \mu,
\end{equation}
\begin{equation} \label{eq:str-const}
c_{ij}^{k} :=\sum_{l=1}^{\mu} \eta^{kl} 
\partial_{i} \partial_{j} \partial_{l} \F, \quad 
(\eta^{ij}) = (\eta_{ij})^{-1}, \quad
i,j,k = 1, \cdots, \mu.
\end{equation}
\subsection{First structure connection of the Frobenius manifold}
For a Frobenius manifold $M$, one can associate 
a connection $\widehat{\ns}$ on $\T_{\PP^{1}\times M}$ \cite{d:1}:
\begin{subequations}
\begin{eqnarray}
\widehat{\ns}_{\delta'} \delta & = & 
\ns_{\delta'} \delta - \frac{1}{u} 
\delta \circ \delta',  \hspace{+1.em}
\delta, \delta' \in \T_{M}, \label{eq:connection-1} \\
\widehat{\ns}_{u\frac{d}{du}} \delta & = &
\frac{1}{u} {\U}(\delta) 
-  {\V}(\delta), \hspace{+1.em} 
\delta \in \T_{M}, \label{eq:connection-2} \\
\widehat{\ns}_{\delta} \frac{d}{du}  & = & 
\widehat{\ns}_{\frac{d}{du}} \frac{d}{du} = 0, \hspace{+1.em}
\delta \in \T_{M}. \label{eq:connection-3}
\end{eqnarray}
\end{subequations}
Here $u$ is the coordinate of $\PP^{1}$, and ${\U}$, ${\V}$ 
are following $\O_{M}$-linear operators acting on ${\T}_{M}$:
\begin{equation} \label{eq:UV}
{\U}(\delta) = E \circ \delta, \quad
{\V}(\delta) = \ns_{\delta} E - \frac{2-d}{2} \delta.
\end{equation} 
\begin{prop}[Proposition 3.1 in \cite{d:1}]\label{prop:first structure connection}
The connection $\widehat{\ns}$ is flat.
\end{prop} 
Note that the parameter $z$ in \cite{d:1} is $1/u$ in this paper. 
The flat connection $\widehat{\ns}$ is called the 
{\it first structure connection} of the Frobenius manifold $M$. 
Let $\varphi = \varphi(t,u)$ be a function on an open subset in $M \times \PP^{1}$. 
We say that $\varphi$ is {\it $\widehat{\ns}$-flat} if it satisfies 
\[
\widehat{\ns} \biggl( \sum_{i=1}^{\mu} 
u\frac{\partial \varphi}{\partial t_{i}} dt_{i} \biggr) = 0
\]
under the identification $\T_{M} \cong \Omega_M^1$ by the non-degenerate 
bilinear form $\eta$. That is, the gradient 
$\Phi = {}^{T} (u\partial_{1}\varphi,\cdots,u\partial_{\mu}\varphi)$
of a $\widehat{\ns}$-flat function $\varphi$ satisfies 
\begin{subequations} \label{eq:qde}
\begin{equation}\label{eq:qde-1}
\frac{\partial}{\partial t_{i}} \Phi = \frac{1}{u} C_{i} \Phi, 
\hspace{+1.em} i = 1, \cdots, \mu,  
\end{equation}
\begin{equation} \label{eq:qde-2}
u\frac{d}{du} \Phi + \frac{1}{u} U \Phi -  V \Phi = 0.
\end{equation}
\end{subequations}
Here $(C_{i})_{jk} = c_{ij}^{k}$, $U$ and $V$ are the matrices representing 
the $\O_{M}$-linear operators $\U$ and $\V$ respectively in the above identification:
\begin{equation} \label{eq:UVmatrix}
{\mathcal U}(dt_{i}) = \sum_{k=1}^{\mu} U_{ki} dt_{k}, \quad
{\mathcal V}(dt_{i}) = \sum_{k=1}^{\mu} V_{ki} dt_{k}, \quad i = 1,\dots,\mu.
\end{equation}
The equation \eqref{eq:qde-2} can be considered as a family of 
meromorphic differential equations on $\PP^{1}$ parametrized by 
points on $M$, and Proposition \ref{prop:first structure connection}
implies that this family is isomonodromic. The equation has a regular singular point 
at $u = \infty$ and an irregular singular point of Poincar\'e rank one at $u = 0$.
\subsection{Semi-simple Frobenius manifold and the canonical coordinate}
In this section we recall the notion of semi-simple Frobenius manifolds. 
\begin{defn}
Let $M$ be a Frobenius manifold.
\begin{enumerate}
\item
A point on $t \in M$ is called {\it semi-simple} if there are no nilpotent 
elements in the product $\circ$ on the tangent space $T_{t}M$. 
\item
A Frobenius manifold is called {\it semi-simple} if general points are semi-simple. 
\item
The set $\K:=\{t\in M~|~ t~\text{is not semi-simple}\}$ is called the set of {\it caustics} 
of the Frobenius manifold $M$.
\end{enumerate}
\end{defn}
The semi-simplicity is an open property, namely, the set $M\setminus \K$ is an open subset of $M$. 
\begin{prop}[Theorem~3.1 in \cite{d:2}]
Near a semi-simple point the eigenvalues $\{w_{a}\}_{a=1}^{\mu}$ of the matrix $U$ 
give local coordinates of $M\setminus \K$. They satisfy the following$:$ 
\begin{subequations}
\begin{eqnarray}
\frac{\partial}{\partial w_{a}} \circ \frac{\partial}{\partial w_{b}} &=& 
\delta_{ab} \frac{\partial}{\partial w_{a}},\quad a,b=1,\dots,\mu, \\
e&=&\sum_{a=1}^{\mu} \frac{\partial}{\partial w_{a}},\\
E&=&\sum_{a=1}^{\mu}w_{a}\frac{\partial}{\partial w_{a}}\label{eq:euler}. 
\hspace{+3.em}
\end{eqnarray}
\end{subequations}
\end{prop}
\begin{defn}
The local coordinates $(w_{1},\dots,w_{\mu})$ of $M\setminus \K$ are called the {\it canonical coordinates} 
of the Frobenius manifold $M$.
\end{defn}
We recall the following important fact$:$
\begin{prop}[Corollary~3.1 in \cite{d:2}]
All the points $t\in M$ where the eigenvalues $\{w_{a}\}_{a=1}^{\mu}$ of the matrix $U$ 
are pairwise distinct are semi-simple. 
\end{prop}
Set $\B:=\{t\in M~|~\text{some of eigenvalues of the matrix } U~\text{coincide}\}$ and call $\B$ the {\it bifurcation set} 
of the Frobenius manifold $M$. 
It follows from the above Proposition that $(M\setminus \B)\subset (M\setminus\K)$.
\subsection{Stokes matrix of the first structure connection}
Fix a point on $M\setminus \B$, where the values of canonical coordinates are 
pairwise distinct. Define the $\mu\times\mu$ matrix $\Psi=(\Psi_{ai})_{a,i=1,\dots,\mu}$ by 
\begin{equation} \label{eq:matrix Psi}
\Psi_{ai} :=\frac{\partial w_a}{\partial t_{i}}\cdot 
\eta\left(\frac{\partial}{\partial w_{a}},\frac{\partial}{\partial w_{a}}
\right)^{\frac{1}{2}},\quad a,i=1,\dots,\mu.
\end{equation}
It follows from the equation \eqref{eq:euler} that the matrix $\Psi$ diagonalizes the matrix $U$:
\begin{equation}
\Psi^{-1} U \Psi = {\rm diag}(w_{1},\dots,w_{\mu}).
\end{equation}
We can construct a formal solution of \eqref{eq:qde-2} near $u=0$ at the point as follows.
\begin{prop}[Lemma 4.3 in \cite{d:2}]
There exists a formal matrix fundamental solution  
of the differential equation \eqref{eq:qde-2} in the form 
\begin{equation} \label{eq:formal sol}
\Phi_{{\rm formal}}(u)=\Psi\hspace{+.1em}G(u)\hspace{+.1em}{\rm exp}({W}/{u}).
\end{equation}
Here $\Psi$ is given by \eqref{eq:matrix Psi}, $W={\rm diag}(w_{1},\dots,w_{\mu})$ 
and $G(u)=1+G_{1}u+G_{2}u^{2}+\cdots$ is a $\mu \times \mu$ matrix-valued 
formal power series satisfying 
\[
G^{T}(-u) G(u) = 1.
\]
Such $G(u)$ is unique.
\end{prop}
Since $u=0$ is an irregular singular point of Poincar\'e rank one, 
the formal power series $G(u)$ is a divergent power series in general.  
\begin{defn} \label{admissible ray}
For $0 \le \phi < \pi$, a line 
$\ell = \{u \in \CC^{\times} | \arg{u} = \phi,\phi-\pi \}$ 
is called {\it admissible} if the line through $w_{a}$ and $w_{b}$ is never 
orthogonal to $\ell$ for any $a,b \in \{1, \cdots, \mu \}$ with $a \ne b$. 
\end{defn}
Fix such a line $\ell$, and chose a small number $\varepsilon > 0$ so that 
any line passing through the origin with angle between  $\phi - \varepsilon$ 
and $\phi + \varepsilon$ is admissible. Define sectors in $u$-plane by 
\begin{eqnarray}
D_{\rm right} & = & \{u \in \CC^{\times}~|~ 
\phi - \pi - \varepsilon < \arg{u} < \phi + \varepsilon \}, \\
D_{\rm left} & = & \{u \in \CC^{\times}~|~
\phi - \varepsilon < \arg{u} < \phi + \pi + \varepsilon \}. 
\end{eqnarray} 
The following statement is a consequence of the general theory 
for ordinary differential equations. 
\begin{prop}[e.g., Theorem A in \cite{BLJ}]
There exists unique solutions $\Phi_{\rm right/left}$ 
of \eqref{eq:qde-2} analytic in $u$ in the sectors 
$D_{\rm right/left}$ having the following asymptotic properties:
\begin{eqnarray}
\Phi_{\rm right}(u) & \sim & \Phi_{\rm formal}(u)\quad
\text{as} \hspace{+.5em} u \rightarrow 0,  \hspace{+.3em} u \in D_{\rm right}, \\
\Phi_{\rm left}(u) & \sim & \Phi_{\rm formal}(u)\quad
\text{as} \hspace{+.5em} u \rightarrow 0, \hspace{+.3em} u \in D_{\rm left}.
\end{eqnarray}
\end{prop}
In the sector 
\begin{equation}
D_{+} = \{u \in \CC^{\times}~|~\phi-\varepsilon<\arg{u}<\phi+\varepsilon \} 
\subset D_{\rm right} \cap D_{\rm left}
\end{equation}
we have two analytic solutions $\Phi_{\rm right}$ and $\Phi_{\rm left}$. 
They must be related as 
\begin{equation}
\Phi_{\rm left}(u) = \Phi_{\rm right}(u) S, \quad u \in D_{+}
\end{equation}
with a matrix $S$ independent of $u$. 
The following proposition is a consequence of the isomonodromic property 
of the differential equation \eqref{eq:qde-2}.
\begin{prop}[Theorem 4.4 in {\cite{d:2}}]
The matrix $S$ is locally constant as a function on $M\setminus\B$.
\end{prop}
\begin{defn}
The matrix $S$ is called the {\it Stokes matrix} of the 
first structure connection of $M$ 
(for the admissible line $\ell$).
\end{defn}
\section{Mirror Isomorphism of Frobenius manifolds}
Let $A:=(a_1,a_2, a_3)$ be a triple of positive integers. Set
\begin{equation}
\mu_A=2+\sum_{k=1}^{3} \left(a_{k}-1\right),
\end{equation}
\begin{equation}
\chi_A:=2+\sum_{k=1}^{3}\left(\frac{1}{a_k}-1\right).
\end{equation}
We shall only consider $A$ satisfying the condition $\chi_A>0$ in this paper. 
\subsection{Gromov--Witten theory for $\PP^1_A$}
Following Geigle--Lenzing (cf. Section~1.1 in \cite{gl:1}), we shall introduce orbifold projective lines. 
First, we prepare some notations.
\begin{defn}
Define a ring $R_{A}$ by 
\begin{subequations}
\begin{equation} 
R_{A}:=\CC[X_1,X_2,X_3]\left/I\right.,
\end{equation}
where  $I$ is an ideal generated by the homogeneous polynomial
\begin{equation}
X_3^{a_3}-X_2^{a_2}+ X_1^{a_1}.
\end{equation}
\end{subequations}
Denote by $L_A$ an abelian group generated by three letters $\vec{X_i}$, 
$i=1,2,3$ defined as the quotient 
\begin{subequations}
\begin{equation}
L_A:=\bigoplus_{i=1}^3\ZZ\vec{X}_i\left/M_A\right. ,
\end{equation}
where  $M_A$ is the subgroup generated by the elements
\begin{equation}
a_i\vec{X}_i-a_j\vec{X}_j,\quad 1\le i<j\le 3.
\end{equation}
\end{subequations}
\end{defn}
We then consider the following quotient stack$:$
\begin{defn}\label{defn:gl}
Define a stack $\PP^1_{A}$ by
\begin{equation}
\PP^1_{A}:=\left[\left({\rm Spec}(R_{A})\setminus\{0\}\right)/{\rm Spec}({\CC L_A})\right],
\end{equation}
which is called the {\it orbifold projective line} of type $A$. 
\end{defn}
An orbifold projective line of type $A$ is a Deligne--Mumford stack 
whose coarse moduli space is a smooth projective line $\PP^1$.
For $g\in\ZZ_{\ge 0}$, $n\in\ZZ_{\ge 0}$ 
and $\beta\in H_2(\PP^{1}_{A},\ZZ)$, 
the moduli space (stack) $\overline{\M}_{g,n}(\PP^1_A,\beta)$ 
of orbifold (twisted) stable maps of genus $g$ 
with $n$-marked points of degree $\beta$ is defined.
There exists a virtual fundamental class $[\overline{\M}_{g,n}(\PP^1_A,\beta)]^{vir}$ 
and Gromov--Witten invariants of genus $g$ 
with $n$-marked points of degree $\beta$ are defined as usual except for 
that we have to use the orbifold cohomology group $H^*_{orb}(\PP^1_A,\QQ)$:
\[
\left<\Delta_1,\dots, \Delta_n\right>_{g,n,\beta}^{\PP^1_A}:=
\int_{[\overline{\M}_{g,n}(\PP^1_A,\beta)]^{vir}}ev_1^*\Delta_1\wedge \dots \wedge 
ev_n^*\Delta_n,\quad \Delta_1,\dots,\Delta_n\in H^*_{orb}(\PP^1_A,\QQ),
\]
where $ev^*_i:H^*_{orb}(\PP^1_A,\QQ)\longrightarrow 
H^*(\overline{\M}_{g,n}(\PP^1_A,\beta),\QQ)$ denotes the induced 
homomorphism by the evaluation map. 
We also consider the generating function 
\[
\F_g^{\PP^1_A}:=\sum_{n,\beta}\frac{1}{n!}
\left<{\bf \Delta},\dots, {\bf \Delta}\right>_{g,n,\beta}^{\PP^1_A},
\quad {\bf \Delta}=\sum_{i=1}^{\mu_A}t_i\Delta_i
\]
and call it the genus $g$ potential where 
$\{\Delta_1=1,\Delta_2,\dots, \Delta_{\mu_A}\}$ denotes a $\QQ$-basis of 
$H^*_{orb}(\PP^1_A,\QQ)$. Note that $\F_0^{\PP^1_A}|_{t_1=0}$ is a polynomial 
in $t_2,\dots,t_{\mu_A-1}, e^{t_{\mu_A}}$ since we have 
the divisor axiom and we assume that $\chi_A>0$.
The Gromov--Witten theory for orbifolds developed by 
Abramovich--Graber--Vistoli \cite{agv:1} and Chen--Ruan \cite{cr:1} 
gives us the following.
\begin{prop}\label{prop:GW to Frob}
The genus zero potential $\F_0^{\PP^1_A}$ satisfies 
the Witten--Dijkgraaf--Verlinde--Verlinde $($WDVV$)$ equations.
In particular, there exists a structure of a Frobenius manifold of 
rank $\mu_A$ and dimension one on $M:=\CC^{\mu_A-1}\times (\CC\setminus\{0\})$  
whose non-degenerate symmetric $\O_M$-bilinear form $\eta$ on $\T_M$ 
is given by the Poincar\'{e} pairing.
\end{prop}
\begin{pf}
See Theorem 6.2.1 in \cite{agv:1} and Theorem 3.4.3 in \cite{cr:1}.
\qed
\end{pf}
For simplicity, we shall denote by $M_{\PP^1_A}$ the complex manifold $M$ 
together with the Frobenius structure on $M$ obtained in 
Proposition~\ref{prop:GW to Frob} and call it {\it the Frobenius manifold 
constructed from the Gromov--Witten theory for $\PP^1_A$}.
\subsection{Theory of primitive forms for $f_A$}
Consider an {\it cusp polynomial} of type $A$, namely, a polynomial 
$f_A({\bf x})\in\CC[x_1,x_2,x_3]$ given as
\begin{equation}
f_A({\bf x}):=x_1^{a_1}+x_2^{a_2}+x_3^{a_3}-q^{-1}\cdot x_1 x_2 x_3
\end{equation}
for some $q\in\CC\setminus\{0\}$.
One can easily show that the $\CC$-vector space
\[
\CC[x_1,x_2,x_3]\left/\Big(\frac{\p f_A}{\p x_1},\frac{\p f_A}{\p x_2},
\frac{\p f_A}{\p x_3}\Big)\right.
\]
is of dimension $\mu_A$ since we assume that $\chi_A>0$.
We can consider the {\it universal unfolding} of $f_A$,  
a deformation $F_A$ of $f_A$ defined on $\CC^{3}\times M$, 
$M:=\CC^{\mu_A-1}\times (\CC\setminus\{0\})$ 
over a $\mu_A$-dimensional parameters $({\bf s},s_{\mu_A})\in M$ given as follows:
\begin{defn}
Define a function $F_A({\bf x};{\bf s}, s_{\mu_A})$ defined on $\CC^{3}\times M$ as follows$;$ 
\begin{equation} \label{eq:FA}
F_A({\bf x};{\bf s},s_{\mu_A}):=x_1^{a_1}+x_2^{a_2}+x_3^{a_3}-s_{\mu_A}^{-1}
\cdot x_1 x_2 x_3+s_1\cdot 1+\sum_{i=1}^{3} \sum_{j=1}^{a_i-1}s_{i,j}\cdot x_i^j.
\end{equation}
\end{defn}
Denote by
\[
p:\CC^{3}\times M\longrightarrow M,\quad
({\bf x};{\bf s},s_{\mu_A})\mapsto ({\bf s},s_{\mu_A})
\]
the projection map from the total space to the deformation space.
Set 
\begin{equation}
p_*\O_\C:=\O_{M}[x_1,x_2,x_3]
\left/\Big(\frac{\p F_A}{\p x_1},\frac{\p F_A}{\p x_2},\frac{\p F_A}{\p x_3}\Big)\right..
\end{equation}
$p_*\O_\C$ can be thought of as the direct image of 
the sheaf of relative {\it algebraic} functions on the relative critical 
set $\C$ of $F_A$ with respect to the projection 
$p:\CC^3\times M\longrightarrow M$.
\begin{prop}[Proposition~2.4 in \cite{ist:2}]\label{unfolding}
The function $F_A({\bf x};{\bf s},s_{\mu_A})$ satisfies the following conditions$:$
\begin{enumerate}
\item $F_A({\bf x};{\bf 0},q)=f_A({\bf x})$.
\item The $\O_M$-homomorphism $\rho$ called the Kodaira--Spencer map defined as 
\begin{equation}\label{universal}
\rho: \T_M\longrightarrow p_*\O_\C,\quad \delta\mapsto \delta F_A,
\end{equation}
is an isomorphism.
\end{enumerate}
\end{prop}
\begin{defn}\label{defn:Jacobi ring}
We shall denote by $\circ$ the induced product structure on $\T_M$ 
by the $\O_M$-isomorphism~\eqref{universal}. 
Namely, for $\delta,\delta'\in \T_M$, we have
\begin{equation}
{(\delta\circ\delta')}F_A = 
{\delta}F_A\cdot{\delta'}F_A\ \text{in}\ p_*\O_\C.
\end{equation}
\end{defn}
\begin{defn}\label{defn:vector fields}
The vector field $e$ and $E$ on $M$ corresponding to the unit 
$1$ and $F$ by the $\O_M$-isomorphism~\eqref{universal} 
is called the {\it primitive vector field} and the {\it Euler vector field}, respectively. 
That is, 
\begin{equation}
{e}F_A=1\ \text{and}\ {E}F_A=F_A\ 
\text{in}\ p_*\O_\C.
\end{equation}
\end{defn}
Note that the primitive vector field $e$ and the Euler vector field $E$ on $M$ are given by 
\begin{equation}
e=\frac{\p}{\p s_1},\quad E=s_1\frac{\p }{\p s_1}+
\sum_{i=1}^3\sum_{j=1}^{a_i-1}\frac{a_i-j}{a_i}s_{i,j}\frac{\p}{s_{i,j}}+
\chi_A s_{\mu_A}\frac{\p}{\p s_{\mu_A}}.
\end{equation}
One can then construct the filtered de Rham cohomology group $\H_{F_A}$ 
(whose increasing filtration is denoted by $\H_{F_A}^{(p)}$ $(p\in\ZZ)$), 
the Gauss--Manin connection $\nabla$ on $\H_{F_A}$ and 
the higher residue pairings $K_{F_A}$ on $\H_{F_A}$, which are necessary 
to define a notion of a primitive form. In this paper, we omit the details 
about these objects and refer the interested reader to \cite{ist:2, st:1} 
A primitive form is obtained by Ishibashi--Shiraishi--Takahashi in \cite{ist:2}:
\begin{thm}[Theorem~3.1 in \cite{ist:2}]\label{thm:primitive form}
The element 
\begin{equation}
\zeta_A:=[s_{\mu_A}^{-1}dx_1 \wedge dx_2 \wedge dx_3]\in \H_{F_A}^{(0)}
\end{equation}
is a primitive form for the tuple $(\H_{F_A}^{(0)},\nabla, K_{F_A})$ 
with the minimal exponent $r=1$. 
\end{thm}
Once we have a primitive form $\zeta_A$, we obtain a Frobenius structure on 
$M$ by the general theory developed by K.~Saito. 
\begin{cor}[Corollary~3.2 in \cite{ist:2}]\label{cor:prim to Frob}
The primitive form $\zeta_A$ determines a Frobenius structure of rank $\mu_A$ 
and dimension one on the deformation space $M$ of the universal unfolding of $f_A$. 
More precisely, the non-degenerate symmetric bilinear form $\eta$ on $\T_M$ defined by 
\begin{equation}
\eta(\delta,\delta'):=u^{-3}K_{F_A}(u\nabla_\delta\zeta_A,u\nabla_{\delta'}\zeta_{A}),
\quad \delta,\delta'\in\T_M,
\end{equation}
together with the product $\circ$ on $\T_M$, the primitive vector field 
$e\in\Gamma(M,\T_M)$ and the Euler vector field $E\in\Gamma(M,\T_M)$ 
define a Frobenius structure on $M$ of rank $\mu_A$ and dimension one.
\end{cor}
For simplicity, we shall denote by $M_{f_A,\zeta_A}$ the deformation space $M$ 
together with the Frobenius structure on $M$ obtained in 
Corollary~\ref{cor:prim to Frob} and call it 
{\it the Frobenius manifold constructed from the pair $(f_A,\zeta_A)$}.
\subsection{Mirror isomorphism}
\begin{thm}[Corollary~4.5 in \cite{ist:2}]\label{thm:Frob mirror}
There exists an isomorphism of Frobenius manifolds between $M_{\PP^{1}_A}$ and $M_{f_A,\zeta_A}$. 
\end{thm}
\begin{rem}
A part of the statement is given by Milanov--Tseng 
\cite{mt:1} for the case $a_1=1$ and by Rossi \cite{r:1} 
for the case $\chi_A>0$. 
\end{rem}
\begin{cor}\label{cor:semi-simplicity}
The quantum cohomology ring of the orbifold $\PP^1_A$ is semi-simple.
\end{cor}
\begin{pf}
The statement easily follows from the fact that 
the Frobenius structure constructed from the pair 
$(f_A,\zeta_A)$ is semi-simple whose canonical coordinates 
are given by the critical values.
\qed
\end{pf}
Let $t=(t_{1},\dots,t_{\mu_{A}})$ be the flat coordinate of $M_{\PP^{1}_{A}}$
such that $E=\displaystyle \sum_{i=1}^{\mu_A} 
E_{i} \frac{\partial}{\partial t_{i}} = \sum_{i=1}^{\mu_A}
(1-q_i)t_i\frac{\p}{\p t_i}+\chi_A\frac{\p}{\p t_{\mu_A}}$ where $q_1=0$.
In the above flat coordinate of $M_{\PP^{1}_{A}}$
the matrices in \eqref{eq:UVmatrix} are given by 
\begin{equation} \label{eq:UVm}
{U}_{ij} = \sum_{k = 1}^{\mu_{A}} 
E_{k} c_{ik}^{j}, \quad
{V}_{ij} = \Bigl(q_{i} - \frac{1}{2}\Bigr) \delta_{ij},
\quad i,j = 1, \cdots, \mu_{A}.
\end{equation}
See Section 2 in \cite{d:2} for detail.
Another important corollary of Theorem~\ref{thm:primitive form} 
is the following$:$
\begin{cor}\label{cor:quantum diff}
Let $\widehat{\ns}$ be the first structure connection of the Frobenius manifold 
$M_{\PP^{1}_{A}}$. On a neighborhood of a semi-simple point of $M_{\PP^{1}_{A}}$, 
the oscillatory integral 
\begin{equation}\label{eq:osillatory}
{\mathcal I}({t},u) = (2\pi u)^{-\frac{3}{2}} 
\int_{\Gamma(u)}e^{{F_A}({\bf x};{\bf t},e^{t_{\mu_{A}}})/{u}}\zeta_A 
\end{equation}
gives a $\widehat{\ns}$-flat function, where ${\bf t}=(t_{1},\dots,t_{\mu_{A}-1})$. Here 
$\Gamma(u) \in H_{3}(\CC^{3};{\rm Re}(F_{A}(\bullet; {\bf t},e^{t_{\mu_{A}}})/u) 
\ll 0; \ZZ)$ is a Lefschetz thimble defined in Section \ref{section:Lefschetz thimble}.
\end{cor}
\begin{pf}
Due to the definition {\rm (P4)} of the primitive form $\zeta_A$ 
in Definition~2.38 in \cite{ist:2} and Lemma~3.4 in \cite{ist:2},
we can show that the oscillatory integral \eqref{eq:osillatory}
satisfies the following system of differential equations:
\begin{subequations}
\begin{equation}\label{P4'}
u\frac{\p}{\p t_{i}}\frac{\p}{\p t_{j}}{\mathcal I}
=\left(\frac{\p}{\p t_{i}}\circ 
\frac{\p}{\p t_{j}}\right){\mathcal I},\quad
i,j=1,\dots, \mu_A.
\end{equation}
\begin{equation}\label{P5'}
\left(u\frac{d}{d u}+E\right)\left(u\frac{\p}{\p t_{i}}{\mathcal I}\right)
=\left(q_i-\frac{1}{2}\right)
\left(u\frac{\p}{\p t_{i}}{\mathcal I}\right),\quad i=1,\dots, \mu_A.
\end{equation}
\end{subequations}
The differential equation~\eqref{P4'} and \eqref{eq:prod} imply that 
the gradient of ${\mathcal I}$ satisfies \eqref{eq:qde-1}. 
Since it follows from \eqref{P4'} that 
\[
E\left(u\frac{\partial}{\partial t_{i}}{\mathcal I}\right) = 
\Bigl(E \circ \frac{\partial}{\partial t_{i}} \Bigr)
{\mathcal I}, \quad i=1,\dots,\mu_{A},
\]
the equation~\eqref{eq:qde-2} follows from \eqref{P5'} and \eqref{eq:UVm}. \qed
\end{pf}

\section{Homological Mirror Symmetry}
\subsection{Derived directed fukaya category $D^b{\rm Fuk}^{\to}(f_A)$ for $f_A$}
We regard $F_A$ as a globally defined tame polynomial on $\CC^3\times M$ and then
consider the derived category $D^b{\rm Fuk}^\to(X_{w;{\bf t},e^{t_{\mu_{A}}}})$ 
of the directed Fukaya category ${\rm Fuk}^\to(X_{w;{\bf t},e^{t_{\mu_{A}}}})$. 
Here, ${\rm Fuk}^{\to}(X_{w;{\bf t},e^{t_{\mu_{A}}}})$ is a directed 
$A_\infty$-category which can be thought of as a ``categorification" of  
a distinguished basis of vanishing cycles in the affine variety 
$X_{w;{\bf t},e^{t_{\mu_{A}}}}$ at a point 
$(w;{\bf t},e^{t_{\mu_{A}}})\in(\CC\times M)\setminus \check{\D}$ defined as 
\begin{equation}
X_{w;{\bf t},e^{t_{\mu_{A}}}}:=
\{{\bf x}\in\CC^3~|~F_A({\bf x};{\bf t},e^{t_{\mu_{A}}})=w\},
\end{equation}
\begin{equation}
\check{\D}:= \{ (w;{\bf t},e^{t_{\mu_{A}}}) \in \CC\times M~|~
X_{w;{\bf t},e^{t_{\mu_{A}}}}~\text{is singular} \}. 
\end{equation}
In this paper, we omit the details about Fukaya categories and refer 
the interested reader to \cite{se:1, se:2} for 
Fukaya categories associated to singularities and to \cite{fooo:1} for generality. 
For the convenience of the reader, we give a rough definition of 
${\rm Fuk}^\to(X_{w;{\bf t},e^{t_{\mu_{A}}}})$:
\begin{defn}
The directed Fukaya category ${\rm Fuk}^\to(X_{w;{\bf t},e^{t_{\mu_{A}}}})$ 
is a strictly unital $A_\infty$-category consists of 
\begin{itemize}
\item 
$\mu_A$ vanishing graded Lagrangian submanifolds 
$\L_1,\dots,\L_{\mu_A}$ in the affine variety $X_{w;{\bf t},e^{t_{\mu_{A}}}}$ together with 
an ordering of these objects as $(\L_1,\dots,\L_{\mu_A})$ such that
\begin{equation}
{\rm Fuk}^\to(X_{w;{\bf t},e^{t_{\mu_{A}}}})(\L_i,\L_j)=
\begin{cases}
0 \quad &\text{if}\quad  i> j,\\
\CC\cdot {\rm id}_{\L_i}\quad &\text{if}\quad  i=j,\\
\displaystyle\bigoplus_{p\in \L_i\cap\L_j}\CC [{\rm deg}(p)]\quad &\text{if}\quad  i<j,
\end{cases}
\end{equation}
where $[-]$ denotes the translation of the complex, ${\deg}(p)$ is defined by 
the gradings ${\rm gr}_{\L_i}:\L_i\longrightarrow \RR$ and 
${\rm gr}_{\L_j}:\L_j\longrightarrow \RR$   
as the largest integer less than or equal to ${\rm gr}(\L_j)|_p-{\rm gr}(\L_i)|_p$, 
\item the $($non-trivial$)$ composition maps 
\begin{multline}
\m_A^n: {\rm Fuk}^\to(X_{w;{\bf t},e^{t_{\mu_{A}}}})(\L_{i_{n-1}},\L_{i_n})\otimes_\CC\dots \otimes_\CC 
{\rm Fuk}^\to(X_{w;{\bf t},e^{t_{\mu_{A}}}})(\L_{i_{1}},\L_{i_2}) \\
\longrightarrow  {\rm Fuk}^\to(X_{w;{\bf t},e^{t_{\mu_{A}}}})(\L_{i_{1}},\L_{i_n})[2-n], \quad i_1<\dots <i_n,
\end{multline}
defined by the ``numbers of pseudo-holomorphic polygons" with boundaries on 
$\L_1,\dots,\L_{\mu_A}$ and corners on intersection points.
\end{itemize}
\end{defn}
The $A_\infty$-category ${\rm Fuk}^{\to}(X_{w;{\bf t},e^{t_{\mu_{A}}}})$ depends on many choices other than 
the initial data $f_A$, especially, on the choice of the point $(w;{\bf t},e^{t_{\mu_{A}}})\in(\CC\times M)\setminus\widehat{\D}$.  
However, it turns out that the derived category $D^b{\rm Fuk}^\to (X_{w;{\bf t},e^{t_{\mu_{A}}}})$ of ${\rm Fuk}^{\to}(X_{w;{\bf t},e^{t_{\mu_{A}}}})$ becomes an invariant of the polynomial $f_A$ as a triangulated category, which we shall usually denote by
$D^b{\rm Fuk}^{\to}(f_A)$ for simplicity.
Note also that the ordered set of objects $(\L_1,\dots, \L_{\mu_A})$ 
forms a full exceptional collection in $D^b{\rm Fuk}^\to (f_A)$ by definition.
\subsection{Mirror equivalence}
The following theorem is proven 
by Auroux-Katzarkov-Orlov~\cite{ako:1}, Seidel~\cite{se:1} and van Straten \cite{str:1} 
for the case $A=(1,p,q)$, $p,q\ge 1$ 
and by the second named author~\cite{t:2} for the cases $A=(2,2,r)$, 
$r\ge 2$, $A=(2,3,3)$, $A=(2,3,4)$ and $A=(2,3,5)$.  
\begin{thm}\label{thm:hms}
There exists a triangulated equivalence
\begin{equation}\label{eq:hms}
\sigma:D^b{\rm Fuk}^\to (f_A)\cong D^b{\rm coh}(\PP^1_A).
\end{equation}
In particular, $(\sigma(\L_1),\dots, \sigma(\L_{\mu_A}))$ forms 
a full exceptional collection in $D^b{\rm coh}(\PP^1_A)$ and 
\begin{equation}
\chi(\sigma(\L_i),\sigma(\L_j))+\chi(\sigma(\L_j),\sigma(\L_i))
=-I(\L_i,\L_j),\quad i,j=1,\dots,\mu_A,
\end{equation}
where $I$ is the intersection form on the middle homology group $H_{2}(X_{w;{\bf t},e^{t_{\mu_{A}}}};\ZZ)$ 
of the affine variety $X_{w;{\bf t},e^{t_{\mu_{A}}}}$.
\end{thm}

\section{Lefschetz thimbles and vanishing cycles}
\label{section:Lefschetz thimble}
\begin{figure}[h]
\begin{center}
\includegraphics[width=100mm]{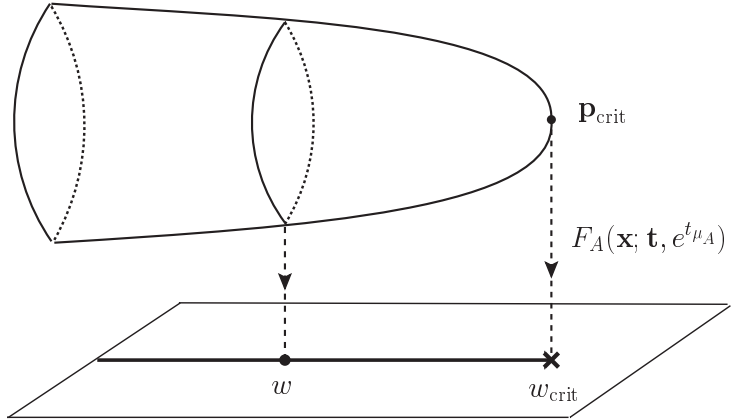}
\end{center} 
\caption{Lefschetz thimble.}
\label{fig:Lefschetz thimble}
\end{figure}

In what follows, we identify $M_{\PP^{1}_{A}}$ and $M_{f_{A},\zeta_{A}}$ 
by the mirror isomorphism of Theorem \ref{thm:Frob mirror}, and denote it by $M$.
Fix a point $({\bf t},e^{t_{\mu_{A}}}) \in M\setminus\B$ and $u \in \CC^{\times}$ 
such that $e^{{\pi\sqrt{-1}}/{2}}u$ lies on an admissible line in the sense of 
Definition \ref{admissible ray}. Then, for each critical point 
${\bf p}_{\rm crit} \in \CC^{3}$ of 
$F_{A}({\bf x};{\bf t},e^{t_{\mu_{A}}})$, we can define a relative $3$-cycle 
\begin{equation}\label{eq:Lef thi}
\Gamma(u) \in H_{3}(\CC^{3},{\rm Re}(F_{A}(\bullet;{\bf t},e^{t_{\mu_{A}}})/u)\ll0;\ZZ),
\end{equation} 
called the {\it Lefschetz thimble} for ${\bf p}_{\rm crit}$ as follows. The image of 
$\Gamma(u)$ by $F_{A}({\bf x}; {\bf t},e^{t_{\mu_{A}}})$ is a half-line starting from 
$w_{\rm crit} = F_{A}({\bf p}_{\rm crit}; {\bf t},e^{t_{\mu_{A}}})$ in the direction of 
$(\pi + \arg u)$, and the fiber above a point $w$ on this half-line is the $2$-cycle 
in $X_{w;{\bf t},e^{t_{\mu_A}}}$ which vanishes at ${\bf p}_{\rm crit}$ by the parallel 
transport along this half-line (see Figure \ref{fig:Lefschetz thimble}). 
Note that all critical points of $F_{A}({\bf x};{\bf t},e^{t_{\,u_{A}}})$ 
are non-degenerate since we take a semi-simple point. 
\begin{prop}\label{iso of homology}
For any point 
$(w;{\bf t},e^{t_{\mu_{A}}})\in\left(\CC\times M\right)\setminus\check {\D}$, 
we have an isomorphism 
\begin{equation}
\p:H_3(\CC^3,X_{w;{\bf t},e^{t_{\mu_{A}}}};\ZZ)\cong 
H_2(X_{w;{\bf t},e^{t_{\mu_{A}}}};\ZZ).
\end{equation}
\end{prop}
\begin{pf}
The relative homology long exact sequence yields the statement. \qed
\end{pf}
For a point $({\bf t},e^{t_{\mu_{A}}}) \in M\setminus\B$ and an admissible $u$, 
take a regular value $w_{0} \in \CC$ with ${\rm Re}(w_{0}/u)$ is small enough so that 
\[
H_{3}(\CC^{3}, {\rm Re}(F_{A}(\bullet; {\bf t},e^{t_{\mu_{A}}})/u)\ll0;\ZZ)
\cong H_3(\CC^3,X_{w_{0};{\bf t},e^{t_{\mu_{A}}}};\ZZ)
\]
holds. Then, the homology class represented by Lefschetz thimbles are 
uniquely characterized by the vanishing cycles in the affine variety 
$X_{w_{0};{\bf t},e^{t_{\mu_{A}}}}$ by Proposition \ref{iso of homology}. 
%

\section{Main Theorem}
\begin{figure}[h]
\begin{center}
\includegraphics[width=70mm]{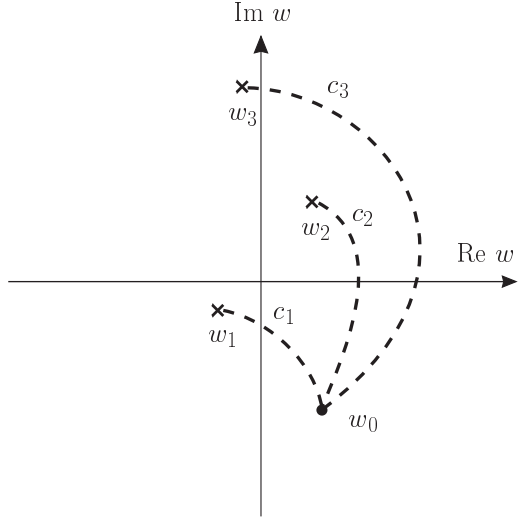}
\end{center} \vspace{-1.em}
\caption{The path $c_{i}$ which determines the vanishing cycle $\L_{i}$.}
\label{fig:distinguished basis}
\end{figure}

In this section we discuss on a neighborhood of a point on $M\setminus\B$. 
Fix an admissible line $\ell = e^{\sqrt{-1}\phi}\RR\setminus\{0\}$ 
($0 \le \phi < \pi$). Let ${\bf p}_{1}, \dots, {\bf p}_{\mu_{A}}$ be the critical points o
f $F_{A}({\bf x}; {\bf t},e^{t_{\mu_{A}}})$, $\Gamma_{1}(u), \cdots, \Gamma_{\mu_{A}}(u)$ 
be the Lefschetz thimbles for these critical points when $e^{\pi\sqrt{-1}/2} u \in \ell$, 
and define 
\begin{equation}\label{eq:osc int}
{\mathcal I}_{i}({t},u) = (2\pi u)^{-\frac{3}{2}} 
\int_{\Gamma_{i}(u)}e^{{F_A}({\bf x};{\bf t},e^{t_{\mu_{A}}})/{u}}\zeta_A,
\quad i = 1,\dots,\mu_{A}. 
\end{equation}
By the saddle-point method, the oscillatory integral has the asymptotic expansion
\begin{equation} \label{saddle-point approximation}
{\mathcal I}_{i}({t},u) = 
\frac{e^{-t_{\mu_{A}}}e^{{w_{i}(t)}/{u}}}
{\sqrt{-\Delta_{i}(t)}}(1 + O(u)),\quad i=1,\dots,\mu_{A},
\end{equation}
as $u \rightarrow 0$ with the fixed argument 
since $s_{\mu_A}=e^{t_{\mu_{A}}}$ 
(see sentences after Lemma~4.3 in \cite{ist:2}). 
Here $w_{i}({t}) = F_{A}({\bf p}_{i};{\bf t},e^{t_{\mu_{A}}})$ 
is the $i$-th critical value, and $\Delta_{i}({t})$ 
is the Hessian at ${\bf p}_{i}$
\begin{equation}
\Delta_{i}(t) := {\rm det}\left( \frac{\partial^{2} F_{A}}
{\partial x_{k} \partial x_{l}}({\bf p}_{i};{\bf t},e^{t_{\mu_{A}}})
\right)_{k,l=1,2,3},\quad i=1,\dots,\mu_{A}.
\end{equation}
Since $\p/\p w_1,\dots, \p/\p w_{\mu_A}$ are basic idempotents, 
it follows that 
\begin{equation}
\frac{\p F_{A}}
{\p w_i}({\bf p}_{i};{\bf t},e^{t_{\mu_{A}}})=1.
\end{equation}
Therefore, by Definition~2.29 and Definition~2.33 (iv) in \cite{ist:2}, we have
\begin{equation}
\eta\left(\frac{\p }{\p w_i},\frac{\p }{\p w_i}\right)=
-\frac{e^{-2t_{\mu_{A}}}}{\Delta_{i}(t)},\quad i=1,\dots,\mu_{A},
\end{equation}
and hence
\begin{equation} \label{eq:psi1}
\Psi_{1i} = \frac{e^{-t_{\mu_{A}}}}{\sqrt{-\Delta_{i}(t)}},\quad i=1,\dots,\mu_{A}.
\end{equation}
The equality \eqref{eq:psi1} and Corollary \ref{cor:quantum diff} imply that 
the asymptotic expansion of the gradient 
$^{T}(u\partial_{1}\I_{i},\dots,u\partial_{\mu_{A}}\I_{i})$ 
coincides with the $i$-th column of the formal matrix solution 
$\Phi_{\rm formal}$ if a suitable branch of the square root is chosen. 
The integration cycles $\Gamma_{i}(u)$ undergoes a discontinuous change 
when $u$ cross a line such that the half-line starting from the critical value 
$w_{i}$ in the direction of $(\pi + \arg u)$ pass through another critical value $w_{j}$.
These discontinuities cause Stokes phenomena for the oscillatory integrals. 
In order to determine the Stokes matrix $S$ of the Frobenius manifold $M$ 
for the admissible line $\ell$, we establish the correspondence between the 
analytic solutions $\Phi_{\rm right/left}$ of \eqref{eq:qde-2} and ${\mathcal I}_{i}$.
\begin{figure}[h]
\begin{center}
\includegraphics[width=80mm]{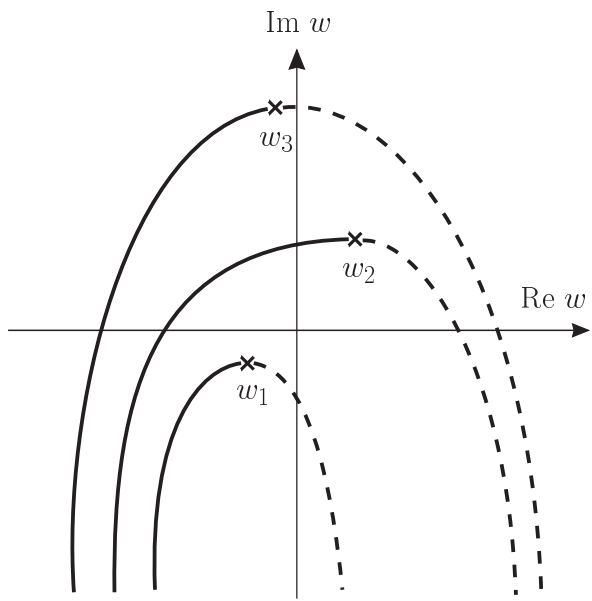}
\end{center}  
\caption{The image of $\Gamma_{i, {\rm right/left}}(u = +\sqrt{-1})$ by $F_{A}$.}
\label{fig:at u=i}
\end{figure}
In what follows, without loss of generality, we assume that 
the positive imaginary axis $\ell = \sqrt{-1}\hspace{+.2em}{\RR}\setminus\{0\}$ 
is admissible, and discuss the Stokes matrix for $\ell$. Take a small 
positive number $\varepsilon$ such that all lines passing through the origin 
with angle between $\pi/2-\varepsilon$ to $\pi/2+\varepsilon$ are admissible.
Order critical values $\{w_{i}\}_{i=1}^{\mu_{A}}$ so that 
\[
|e^{{w_{1}}/{u}}| \ll |e^{{w_{2}}/{u}}| \ll \cdots \ll |e^{{w_{\mu_{A}}}/{u}}|
\]
holds as $u \rightarrow 0$ along the line $\arg u = \pi/2$. Consider the local system 
on $\CC^{\times}$ whose fiber on $u \in \CC^{\times}$ is the relative homology group 
$H_{3}(\CC^{3},{\rm Re}(F_{A}(\bullet;{\bf t},e^{t_{\mu_{A}}})/u) \ll 0;\ZZ)$, 
and let $\Gamma_{i,\rm right/left}$ be a section of the local system on 
$D_{\rm right/left}$ satisfying the following condition; 
\begin{itemize}
\item[(A)] for $u \in \{ u \in \CC^{\times}~|-\varepsilon<\arg u<\varepsilon \} 
\subset D_{\rm right}$ (resp., $u \in \{ u \in \CC^{\times}~|~ 
\pi-\varepsilon<\arg u<\pi+\varepsilon \} \subset D_{\rm left}$),
$\Gamma_{i, {\rm right}}(u)$ (resp., $\Gamma_{i, {\rm left}}(u)$) 
coincides with the relative homology class represented by 
the Lefschetz thimble for the $i$-th critical point. 
\end{itemize}
Since $D_{\rm right/left}$ is simply-connected, the condition (A)
determines the section $\Gamma_{i,\rm right/left}$ uniquely. 
Figure \ref{fig:at u=i} describes the projection of cycles in 
$H_{3}(\CC^{3},{\rm Re}(F_{A}(\bullet;{\bf t},e^{t_{\mu_{A}}})/u) \ll 0;\ZZ)$ 
representing these homology classes when $u=+\sqrt{-1} \in D_{+}$ in the case 
$\mu_{A} = 3$. Solid curves (resp., dotted curves) in Figure \ref{fig:at u=i} 
express the image of cycles representing $\Gamma_{i, {\rm right}}(u=+\sqrt{-1})$ 
(resp., $\Gamma_{i, {\rm left}}(u=+\sqrt{-1})$) by $F_{A}({\bf x};{\bf t},e^{t_{\mu_{A}}})$.
Define 
\begin{eqnarray}
{\mathcal I}_{i,\rm right}(t,u)& = &(2\pi u)^{-\frac{3}{2}}\int_{\Gamma_{i, \rm right}(u)}
e^{{F_{A}({\bf x};{\bf t},e^{t_{\mu_{A}}})}/{u}}\zeta_{A}, \\
{\mathcal I}_{i,\rm left}(t,u)& = &(2\pi u)^{-\frac{3}{2}}\int_{\Gamma_{i, \rm left}(u)}
e^{{F_{A}({\bf x};{\bf t},e^{t_{\mu_{A}}})}/{u}}\zeta_{A}, 
\end{eqnarray}
for $i = 1,\cdots,\mu_{A}$. The matrix valued fuction
\begin{equation}\label{eq:solution matrix}
\left( u\frac{\partial \I_{j, {\rm right/left}}}{\partial t_{i}}(t,u)
\right)_{i,j=1,\dots,\mu_{A}}
\end{equation}
is a fundamental solution of \eqref{eq:qde} which is asymptotic to 
$\Phi_{\rm formal}$ as $u \rightarrow 0$ in the whole sector $D_{\rm right/left}$ 
by the condition (A). Therefore, \eqref{eq:solution matrix} coincides with 
$\Phi_{\rm right/left}$ and the desired Stokes matrix $S$ can be read off 
from the monodromy of integration cycles since the integrand is single-valued. 
That is, for $u \in D_{+}$, $S = (S_{ij})_{i,j=1,\cdots,\mu_{A}}$ satisfies 
\begin{equation}
\Gamma_{j, \rm left}(u) = \sum_{i = 1}^{\mu_{A}} S_{ij} \Gamma_{i,\rm right}(u)
\quad j = 1,\cdots,\mu_{A}.
\end{equation} 
As we discussed in Section \ref{section:Lefschetz thimble}, for a fixed $u \in D_{+}$,
we can take a regular value $w_{0} \in {\CC}$ of $F_{A}({\bf x};{\bf t},e^{t_{\mu_{A}}})$ 
with ${\rm Re}(w_{0}/u)<0$ is small enough such that 
\begin{equation} \label{eq:LT and VC}
H_{3}(\CC^{3},{\rm Re}(F_{A}(\bullet;{\bf t},e^{t_{\mu_{A}}})/u)\ll0;\ZZ) 
\cong H_2(X_{w_{0};{\bf t},e^{t_{\mu_{A}}}};\ZZ)
\end{equation}
holds. Fix such a point $w_{0}$ and take paths $c_{i}$ 
between $w_{0}$ and $w_{i}$ for $i = 1,\dots,\mu_{A}$ satisfying
\begin{itemize}
\item for $i \ne j$, $c_{i}$ and $c_{j}$ have a unique common point $w_{0}$,
\item $c_{1},\dots,c_{\mu_{A}}$ are ordered such that the following holds; 
for $i<j$, $w \in c_{i}$ and $w' \in c_{j}$, ${\rm Re}~w < {\rm Re}~w'$ holds
if ${\rm Im}~w = {\rm Im}~w'$.
\end{itemize}
See Figure \ref{fig:distinguished basis} for an example of such paths.
Denote by $\L_{i}$ the cycle in ${X}_{w_{0};{\bf t},e^{t_{\mu_{A}}}}$ which vanishes at 
$w_{i}$ by the parallel transport along the path $c_{i}$ ($i = 1,\dots,\mu_{A}$). 
Then, the ordered set $({\mathcal L}_{1},\dots,{\mathcal L}_{\mu_{A}})$ of cycles 
form a {\it distinguished basis of vanishing cycles} in the sense in \cite{AGV}. 
%

%
According to the Picard-Lefschetz theory, the relationship between the cycles 
$\{\Gamma_{i,{\rm right}}\}_{i=1}^{\mu_{A}}$ and
$\{\Gamma_{i,{\rm left}}\}_{i=1}^{\mu_{A}}$ are expressed by 
the intersection numbers of these vanishing cycles.
Here we recall the Picard-Lefschetz formula.
\begin{prop}[e.g., Section 2 in \cite{AGV}]\label{prop:PLF}
For $i=1,\dots,\mu_{A}$, let 
$h_{i} \in {\rm Aut}(H_{2}(X_{w_{0};{\bf t},e^{t_{\mu_{A}}}};\ZZ))$ 
be the monodromy operator along the loop 
$\tau_{i} \in \pi_{1}(\CC\setminus\{w_{1},\dots,w_{\mu_{A}}\},w_{0})$, 
which goes along the path $c_{i}$ from $w_{0}$ to $w_{i}$, turns around $w_{i}$ 
in the positive direction (anti-clockwise) and returns to $w_{0}$ along $c_{i}$. 
For any cycle $\L \in H_{2}(X_{w_{0};{\bf t},e^{t_{\mu_{A}}}};\ZZ)$, we have 
\begin{equation} \label{eq:PLF}
h_{i}(\L) = \L + I(\L,\L_{i}) \L_{i}, \quad i=1,\dots,\mu_{A}.
\end{equation}
\end{prop}

\begin{figure}[h]
\begin{minipage}{0.450\hsize} 
\begin{center}
\hspace{-3.em}
\includegraphics[width=80mm]{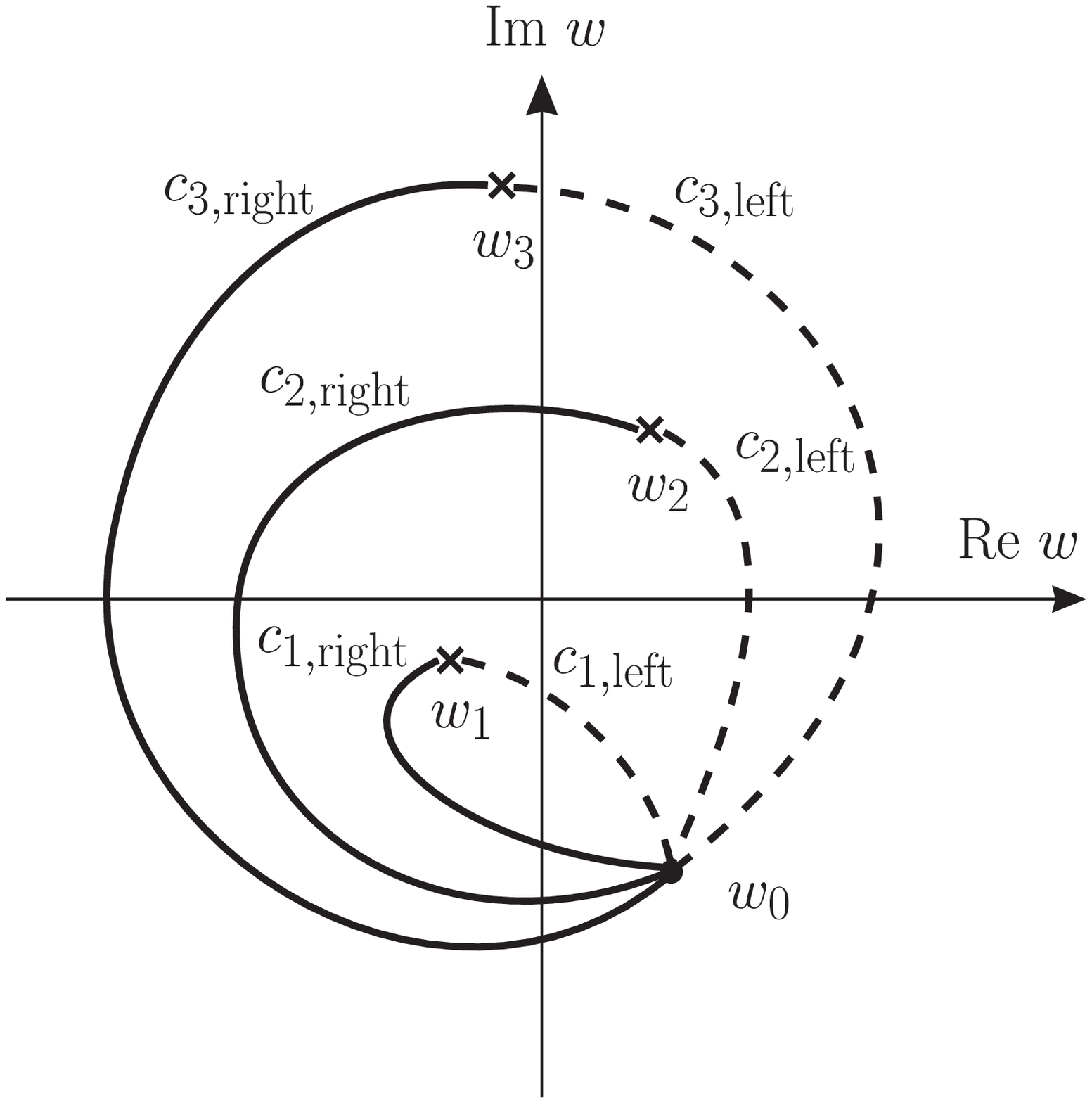}
\end{center}  \vspace{-1.em}
\caption{$c_{i, {\rm right/left}}$.}
\label{fig:Picard-Lefschetz1}
\end{minipage}  \hspace{-.5em}
\begin{minipage}{0.450\hsize}
\begin{center}
\includegraphics[width=80mm]{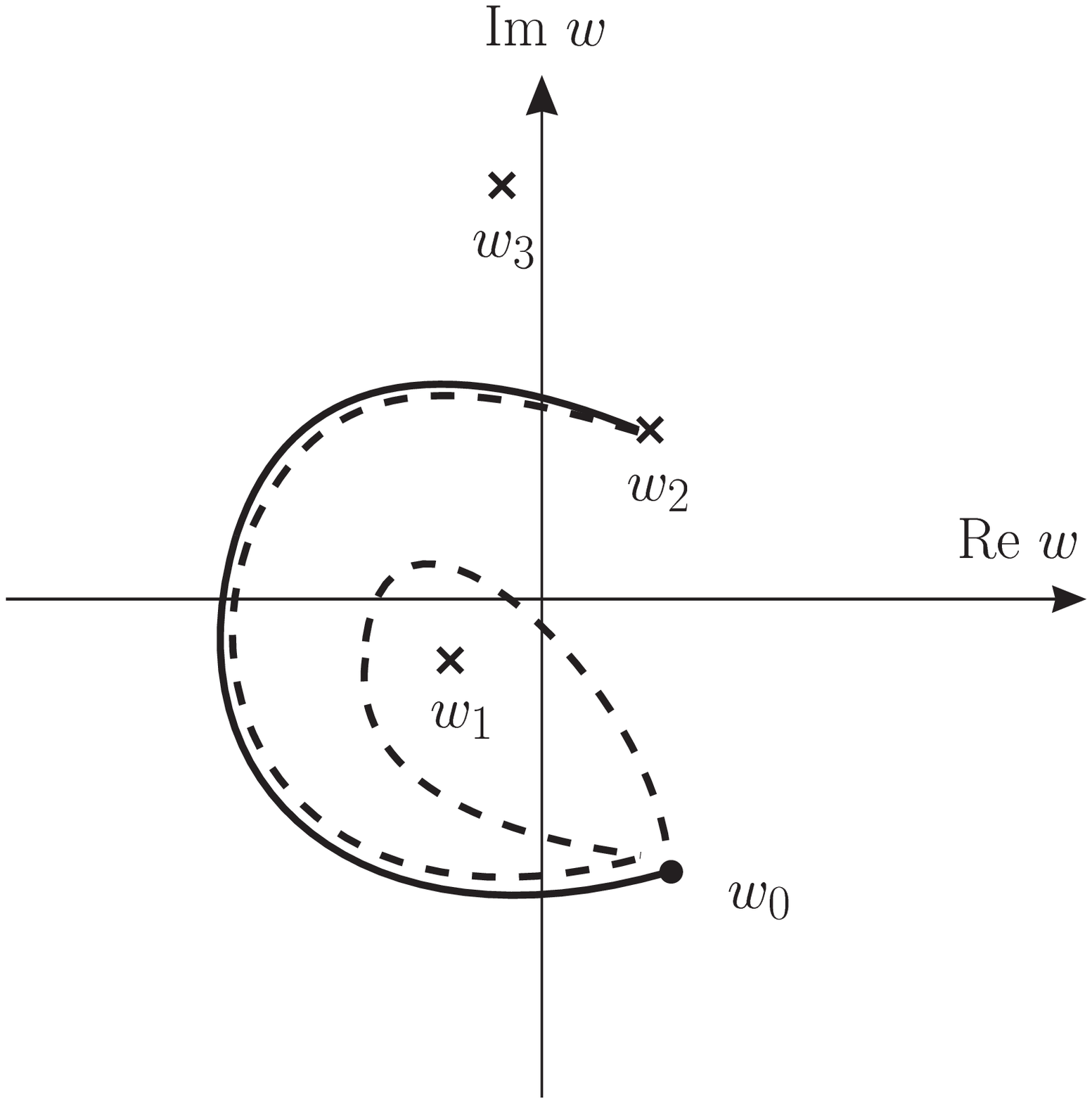}
\end{center}  \vspace{-.4em}
\caption{Deformation of $c_{2, {\rm left}}$.}
\label{fig:Picard-Lefschetz2}
\end{minipage}
\end{figure}

Using the Picard-Lefshetz formula, we obtain the following.
\begin{prop} \label{prop:Picard-Lefschetz}
The following equality holds:
\begin{equation}\label{eq:P-L}
\Gamma_{j,{\rm left}}(u) = \Gamma_{j,{\rm right}}(u) - \sum_{i=1}^{j-1}
I({\mathcal L}_{i},{\mathcal L}_{j}) \Gamma_{i,{\rm right}}(u),\quad j=1,\cdots,\mu_{A}
\end{equation}
for $u \in D_{+}$. Here $I$ is the intersection form on 
$H_{2}(X_{w_{0};{\bf t},e^{t_{\mu_{A}}}};\ZZ)$.
\end{prop}
\begin{pf}
For $u \in D_{+}$ and $i=1,\dots,\mu_{A}$, let 
$\L_{i, {\rm right/left}} \in H_{2}(X_{w_{0};{\bf t},e^{t_{\mu_{A}}}};\ZZ)$ 
be the vanishing cycles corresponding to the homology class 
$\Gamma_{i, {\rm right/left}}(u)$ by the isomorphism \eqref{eq:LT and VC}. 
That is, $\L_{i, {\rm right/left}}$ is the cycle which vanishes at $w_{i}$ 
along the path $c_{i, {\rm right/left}}$ in Figure \ref{fig:Picard-Lefschetz1}, 
and hence $\L_{i} = \L_{i,{\rm left}}$ since $c_{i}$ and $c_{i, {\rm left}}$ 
are homotopic as paths on $\CC\setminus\{w_{1},\dots,w_{\mu_{A}}\}$. 
In view of Figure \ref{fig:Picard-Lefschetz1},  
\begin{equation} \label{eq:L1}
\L_{1,{\rm left}} = \L_{1,{\rm right}}
\end{equation}
holds obviously. For $\L_{2, {\rm right/left}}$, since we can deform the 
path $c_{2,{\rm left}}$ cycles homotopically as in Figure \ref{fig:Picard-Lefschetz2}, 
we have $\L_{2, {\rm right}} = h_{1} (\L_{2, {\rm left}})$. 
It follows from \eqref{eq:PLF} and \eqref{eq:L1} that
\begin{equation}
\L_{2, {\rm left}} =  \L_{2, {\rm right}} - I(\L_{1}, \L_{2}) \L_{1, {\rm right}}.
\end{equation}
Similarly, $\L_{i,{\rm right}}=(h_{1}\circ h_{2}\circ\dots\circ h_{i-1})(\L_{i,{\rm left}})$ 
holds for general $i$. Using the Picard-Lefschetz formula \eqref{eq:PLF} iteratively, 
we obtain 
\begin{equation} \label{eq:PL}
\L_{j, {\rm left}} = \L_{j, {\rm right}} - \sum_{i=1}^{j-1}
I({\mathcal L}_{i},{\mathcal L}_{j}) \L_{j,{\rm right}},\quad j=1,\cdots,\mu_{A}.
\end{equation}
The equality \eqref{eq:P-L} follows from \eqref{eq:PL} and 
the isomorphism \eqref{eq:LT and VC}. \qed
\end{pf}
%

%
By Proposition \ref{prop:Picard-Lefschetz} we have 
\begin{equation}
S_{ij}=
\begin{cases}
0 \quad &\text{if}\quad  i>j,\\
1 \quad &\text{if}\quad  i=j,\\
-I(\L_{i},\L_{j}) \quad &\text{if}\quad  i<j.
\end{cases}
\end{equation} 
That is, the Stokes matrix $S$ and the intersection matrix
$I = (I(\L_{i}, \L_{j}))_{i,j=1,\dots,\mu_{A}}$ are related as 
\begin{equation}\label{Stokes is I}
S + {}^{T}S = -I.
\end{equation}
On the other hand, as a consequence of the homological mirror symmetry 
of Section~4, the intersection matrix $I$ and the Euler matrix 
$\chi = (\sigma(\L_{i}),\sigma(\L_{j}))_{i,j=1,\dots,\mu_{A}}$ 
for the full exceptional collection $(\sigma(\L_{1}),\dots,\sigma(\L_{\mu_{A}}))$
of $D^{b}{\rm coh}(\PP^{1}_{A})$ also satisfy 
\begin{equation} \label{eq:Euler is intersection}
\chi + {}^{T}\chi = -I.
\end{equation}
Since both the Stokes matrix $S$ and the Euler matrix $\chi$ 
are upper triangular matrices with all diagonal entries one,
the relation \eqref{Stokes is I} and \eqref{eq:Euler is intersection} 
implies that $S = \chi$. That is, we have the following statement, 
which shows the Dubrovin's conjecture (\cite{d:3})
for $\PP^{1}_{A}$ with $\chi_{A} > 0$.
\begin{thm}
For any point on $M\setminus\B$ and any admissible line, there exists a f
ull exceptional collection $(\E_{1},\cdots,\E_{\mu_{A}})$ of 
$D^{b}{\rm coh}(\PP^{1}_{A})$ such that the Stokes matrix 
$S = (S_{ij})_{i,j=1,\dots,\mu_{A}}$ of the first structure connection of the 
Frobenius manifold $M_{\PP^1_A}$ coincides with the Euler matrix of them: 
\begin{equation}
S_{ij}={\chi}([\E_{i}],[\E_{j}]), \quad
i,j=1,\dots,\mu_{A}.
\end{equation}
\end{thm}



\begin{thebibliography}{99}
{\small 
\bibitem{agv:1}
D.~Abramovich, T.~Graber and A.~Vistoli,
	{\it Gromov--Witten theory of Deligne--Muford stacks},
	Amer. J. Math. {\bf 130} (2008), no. 5, 1337--1398.
\bibitem{AGV}
V.I.~Arnold, S.M.~Gusein-Zade and A.N.~Varchenko, 
{\it Singularities of Differentiable Maps. Vol.II, 
Monodromy and Asymptotics of Integrals}, 
volume 83 of {\it Monographs in Mathematics}, 
Birkh\"auser Inc., Boston, MA, 1988.  
\bibitem{ako:1}
D.~Auroux, L.~Katzarkov and D.~Orlov,
	{\it Mirror symmetry of weighted projective planes and 
	their noncommutative deformations},
	Annals of Mathematics, {\bf 167} (2008), 867 - 943.
\bibitem{BLJ}
W.~Balser, W.B.~Jurkat and D.A.~Luts, 
{\it Birkhoff invariants and Stokes multipliers for 
meromorphic linear differential equations}, 
J. Math. Anal. Appl. {\bf 71} (1979), 48-94.
\bibitem{cr:1}
W.~Chen and Y.~Ruan,
	{\it Orbifold Gromov--Witten Theory},
	Orbifolds in mathematics and physics (Madison, WI, 2001), 25--85, Contemp. Math., 
	{\bf 310}, Amer. Math. Soc., Providence, RI, 2002.
\bibitem{d:1}
B.~Dubrovin, {\it Geometry of 2d topological field theories},
Integrable systems and quantum groups
(Montecatini Terme, 1993), Lecture Notes in Math., vol. 1620, Springer, 
Berlin, 1996, pp. 120--348.
\bibitem{d:3}
B.~Dubrovin, {\it Geometry and analytic theory of Frobenius manifolds}, In Proceedings 
of the International Congress of Mathematicians, Vol. II (Berlin, 1998), 
number Extra Vol. II, pages 315-326 (electronic), 1998, arXiv:9807034.
\bibitem{d:2}
B.~Dubrovin, {\it Painlev\'e transcendents in two-dimensional topological field theory},  
In ``{\it The Painlev\'e Property: One Century Later}'', 
CRM Ser. Math. Phys., pp 287-412. Springer, New York, 1999.
\bibitem{fooo:1}
	K.~Fukaya, Y.~Oh, H.~Ohta and K.~Ono,
	{\it Lagrangian Intersection Floer Theory: Anomaly and Obstruction}, 
	AMS/IP Studies in Advanced Mathematics, vol 46 I \& II, AMS/International Press, 2009.
\bibitem{gl:1}
       W.~Geigle and H.~Lenzing,
       {\it A class of weighted projective curves arising in
       representation theory of finite-dimensional algebras},
       Singularities, representation of algebras, and vector bundles 
       (Lambrecht, 1985), 9--34, 
        Lecture Notes in Math., 1273, Springer, Berlin, 1987. 
\bibitem{guz}
D.~Guzzetti, {\it Stokes matrices and monodromy groups of the quantum cohomology 
of projective spaces}, Comm. Math. Phys., {\bf 207} (1999), 341-383.
\bibitem{ist:1}
Y. Ishibashi, Y. Shiraishi, A. Takahashi,
{\it A Uniqueness Theorem for Frobenius Manifolds and Gromov--Witten Theory 
for Orbifold Projective Lines}, 
to appear in J. Reine Angew. Math..
\bibitem{ist:2}
Y. Ishibashi, Y. Shiraishi, A. Takahashi,
{\it Primitive Forms for Affine Cusp Polynomials}, 
arXiv:1211.1128.
\bibitem{mt:1}
Todor E. Milanov, Hsian-Hua Tseng,
{\it The spaces of Laurent polynomials, $\mathbb{P}^1$-orbifolds, 
and integrable hierarchies}, Journal f\"{u}r die reine und angewandte 
Mathematik (Crelle's Journal), Volume 2008, Issue 622, Pages189--235.
\bibitem{r:1}
P.~Rossi,
{\it Gromov-Witten theory of orbicurves, the space of tri-polynomials 
and Symplectic Field Theory of Seifert fibrations},
Math. Ann., {\bf 348} (2010), 265--287.
\bibitem{st:1}
K.~Saito and A.~Takahashi,
{\it From Primitive Forms to Frobenius manifolds}, 
Proceedings of Symnposia in Pure Mathematics, {\bf 78} (2008) 31--48.
\bibitem{se:1}
	P.~Seidel,
	{\it More about vanishing cycles and mutation, Symplectic Geometry 
	and Mirror Symmetry}, 
	Proceedings of the 4th KIAS Annual International Conference, 
	World Scientific, 2001, 429-465.
\bibitem{se:2}
	P.~Seidel,
	{\it Fukaya categories and Picard-Lefschetz theory},
	Zurich Lectures in Advanced Mathematics, European Mathematical Society, 2008.
\bibitem{str:1}
	D.~van Straten,
	{\it Mirror Symmetry for $\PP^1$-orbifolds},
	unpublished paper based on talks given in Trieste, Marienthal 
	and G\"oteborg in September 2002.
\bibitem{t:1}
A.~Takahashi, 
{\it Weighted projective lines associated to regular systems of weights of dual type},
Adv. Stud. Pure Math. {\bf 59} (2010), 371--388.
\bibitem{t:2}
A.~Takahashi,
{\it Mirror symmetry between orbifold projective lines and cusp singularities}, 
to appear in Adv. Stud. Pure Math.
\bibitem{u:1}
K.~Ueda, {\it Stokes matrices for the quantum cohomologies of Grassmannians}. 
Int. Math. Res. Not. {\bf 34} (2005), 2075-2086.
\bibitem{u:2}
K.~Ueda, {\it Stokes matrix for the quantum cohomology of cubic surfaces}, 
arXiv:0505350. 
\bibitem{z}
E.~Zaslow, {\it Solitons and Helices: The Search for a Math-Physics Bridge}, 
Comm. Math. Phys, {\bf 175} (1996) 337-375.
}
\end{thebibliography}
\end{document}